\documentclass[english,letterpaper, 10 pt, journal, twoside]{IEEEtran}
\usepackage[T1]{fontenc}
\usepackage[latin9]{inputenc}
\usepackage{amsmath}
\usepackage{amssymb}
\usepackage{graphicx}

\makeatletter
\IEEEoverridecommandlockouts                              

\title{\LARGE \bf
Time-Optimal Collaborative Guidance Using the Generalized Hopf Formula
}

\author{Matthew R. Kirchner, Robert Mar, Gary Hewer, J\'{e}r\^{o}me Darbon, Stanley Osher, Y.T. Chow 
\thanks{This research was supported by the Office of Naval Research, ILIRs 4764 and 5100, and Office of Naval Research grants N00014-16-12119 and N00014-16-12157. }
\thanks{M. Kirchner and G. Hewer are with the Image and Signal Processing Branch, Research Directorate, Code 4F0000D, Naval Air Warfare Center Weapons Division, China Lake, CA 93555, USA  {\tt\small \{matthew.kirchner, gary.hewer\}@navy.mil}}%
\thanks{R. Mar is with the Guidance, Navigation, and Control Branch, Weapons and Energetics Department, Code 472100D, Naval Air Warfare Center Weapons Division, China Lake, CA 93555, USA  {\tt\small robert.t.mar@navy.mil}}%
\thanks{J. Darbon is with the Division of Applied Mathematics, Brown University,
        Providence, RI 02912, USA
        {\tt\small jerome\_darbon@brown.edu}}%
\thanks{S. Osher and Y. T. Chow are with the Department of Mathematics, University of California,
        Los Angeles, CA 90095, USA
        {\tt\small \{sjo, ytchow\}@math.ucla.edu}}%
}

\makeatother

\usepackage{babel}
\begin{document}
\maketitle
\thispagestyle{empty}
\pagestyle{empty} 
\begin{abstract}
Presented is a new method for calculating the time-optimal guidance
control for a multiple vehicle pursuit-evasion system. A joint differential
game of $k$ pursuing vehicles relative to the evader is constructed,
and a Hamilton\textendash Jacobi\textendash Isaacs (HJI) equation
that describes the evolution of the value function is formulated.
The value function is built such that the terminal cost is the squared
distance from the boundary of the terminal surface. Additionally,
all vehicles are assumed to have bounded controls. Typically, a joint
state space constructed in this way would have too large a dimension
to be solved with existing grid-based approaches. The value function
is computed efficiently in high-dimensional space, without a discrete
grid, using the generalized Hopf formula. The optimal time-to-reach
is iteratively solved, and the optimal control is inferred from the
gradient of the value function. 
\end{abstract}

\section{Introduction}

One of the first successful implementations of control laws for pursuit
problems is proportional navigation (PN) \cite{palumbo2010basic},
which attempts to drive the rate of the line-of-sight vector between
pursuer and evading target vehicle to zero. In this derivation, the
target vehicle is assumed moving, but not maneuvering (turning). Generalizations
of this concept attempt to estimate the vehicle maneuver \cite{palumbo2010modern},
but these methods are not optimal since evasion strategy is not considered,
i.e. not formulated as a differential game \cite{isaacs1999differential}.
Additionally, this family of control laws does not account for control
saturation. PN typically requires the magnitude of the control bound
of the pursuer to be much greater than that of the evader to be successful,
on the order of 3-5 times greater \cite{palumbo2010modern}. These
guidance laws are strictly one-on-one in nature, and do not readily
generalize to collaborative systems of multiple vehicles where the
desired pursuit guidance is to 'team' together to capture a target.
These early pursuit problems typically referred to controller designs
as \emph{guidance laws,} and in this letter we will use the terms
controller and guidance interchangeably.

More recently, \cite{pan2012pursuit} proposed a solution to multi-vehicle
pursuit evasion in a plane. In this case the problem was solved sub-optimally
with heuristics in an effort to avoid the computational burden of
direct solution to the Hamilton-Jacobi equation. Additionally, the
method was based on simplified, single-integrator dynamics that require
the vehicles to maneuver instantaneously to ensure capture.

A general alternative is to formulate the pursuit-evasion problem
as a differential game, and derive a Hamilton\textendash Jacobi\textendash Isaacs
(HJI) equation representing the optimal cost-to-go of the system.
Traditionally, numerical solutions to HJI equations require a dense,
discrete grid of the solution space \cite{osher2006level,mitchell2008flexible,mitchell2005time}.
Computing the elements of this grid scales poorly with dimension and
has limited use for problems with dimension of greater than four.
The exponential dimensional scaling in optimization is sometimes referred
to as the ``curse of dimensionality'' \cite{bellman2015adaptive,bellman1957dynamic}.
This phenomenon is seen clearly in \cite{huang2011differential},
which formulated a differential game for a capture-the-flag problem
and solved numerically on a four dimensional grid with \cite{mitchell2004toolbox}.
The computational time was as much as 4 minutes, too slow for real-time
application, even with a coarsely sampled grid of 30 points in each
dimension and with low numeric accuracy. When the grid is increased
to 45 points in each dimension and with high numeric accuracy, the
computation time jumps to an hour.

Recent research \cite{darbon2016algorithms} has discovered numerical
solutions based on the generalized Hopf formula that do not require
a grid and can be used to efficiently compute solutions to a certain
class of Hamilton\textendash Jacobi equations that arise in linear
control theory and differential games. This readily allows the generalization
with pursuit-evasion to collaborative guidance of multiple pursuing
vehicles.

This letter presents a new method for multi-vehicle collaborative
pursuit guidance of a maneuvering target, showing that teams of vehicles
can intercept the target without requiring drastically higher control
bound as in the family of methods in \cite{palumbo2010modern}. A
joint system state space representing the kinematics of all pursuing
vehicles relative to the target was constructed, the dimension of
which makes it infeasible for traditional grid-based methods. This
high-dimensional problem was then efficiently solved using the generalized
Hopf formula, and included the constraint of time-varying bounds on
the magnitude of available vehicle control, while ensuring intercept
when starting within the reachable set.

The rest of the paper is organized as follows. We derive the models
used in the study in Sec. \ref{sec:Pursuit-Evasion-Model} followed
the presentation of efficient solution techniques that employ the
generalized Hopf formula to solve the Hamilton\textendash Jacobi equations
for optimal control and differential games in Sec. \ref{sec:Hamiltonian-Jacobi-Equations}.
The application of these methods to collaborative guidance is given
in Sec. \ref{sec:Collaborative-Guidance}, followed by results on
a planar, multiple vehicle pursuit-evasion game in Sec. \ref{sec:Results}. 

\section{Pursuit-Evasion Model\label{sec:Pursuit-Evasion-Model}}

\subsection{Single Vehicle Model}

\begin{figure}
\begin{centering}
\includegraphics[width=8cm]{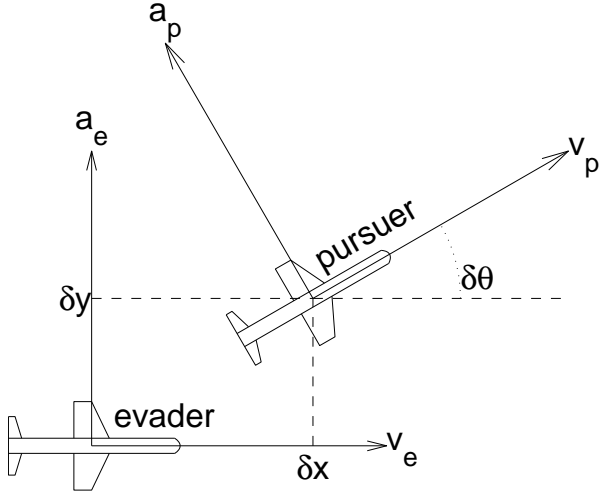}
\par\end{centering}
\caption{The engagement geometry of the system presented in $\left(\ref{eq:Original non-linear system}\right)$.\label{fig:The-engagement-geometry}}
\end{figure}

First consider the pursuit-evasion game with only a single pursuer.
We construct a state space representation of the position and orientation
of the pursuer relative to the evader, with geometry shown in Fig.
\ref{fig:The-engagement-geometry}. With $x=\left[\delta x,\delta y,\delta\theta\right]^{\dagger}$,
the relative system becomes
\begin{equation}
\dot{x}\left(t\right)=\left[\begin{array}{c}
V_{p}\text{cos}\left(\delta\theta\right)-V_{e}+\frac{\delta ya_{e}}{V_{e}}\\
V_{p}\text{sin}\left(\delta\theta\right)-\frac{\delta xa_{e}}{V_{e}}\\
\frac{a_{p}}{V_{p}}-\frac{a_{e}}{V_{e}}
\end{array}\right],\label{eq:Original non-linear system}
\end{equation}
with $V_{p}$ and $V_{e}$ representing the forward speed of the pursuer
and evader, respectively. The terms $a_{p}$ and $a_{e}$ are the
lateral acceleration inputs to the system for each vehicle. These
accelerations are limited by the physical maneuver capabilities of
the vehicles. This system is based on the relative state \cite{mitchell2005time}
of two modified Dubin's car \cite{stipanovic2004decentralized,dubins1957curves}
models, with acceleration instead of the more common turn rate input.
Additionally, we constructed this system to be evader centric, allowing
for the addition of multiple pursuers. Denoting by $\dagger$ the
transpose of a matrix, we introduce the new state vector $x=\left[\delta x,\delta y,\delta v_{x},\delta v_{y}\right]^{\dagger}$,
where $\delta x$ and $\delta y$ are the positional displacement
separating the vehicles (see Figure \ref{fig:The-engagement-geometry}),
$\delta v_{x}=V_{p}-V_{e}$, and $\delta v_{y}$ is the relative vertical
velocity. We proceed to linearize the system $\left(\ref{eq:Original non-linear system}\right)$
with
\begin{align}
\dot{x}\left(t\right) & =\left[\begin{array}{cc}
0_{2} & I_{2}\\
0_{2} & 0_{2}
\end{array}\right]x\left(t\right)+\left[\begin{array}{c}
0\\
0\\
0\\
\pm1
\end{array}\right]a_{p}+\left[\begin{array}{c}
0\\
0\\
0\\
-1
\end{array}\right]a_{e},\label{eq:linearized system}\\
 & =Ax\left(t\right)+Ba_{p}+Da_{e},\nonumber 
\end{align}
with the $\pm$ sign needed depending on whether its tail-chase $\left(+\right)$
or head-on $\left(-\right)$ engagement. The linearization at first
glance may seem extreme, but this linearization strategy is used when
deriving proportional navigation, or its variants such as augmented
proportional guidance and extended proportional guidance, using linear
quadratic control techniques \cite{palumbo2010modern}. The controls
for the pursuer are constrained to the set $\mathcal{A}_{p}=\left\{ a_{p}:\left\Vert Q_{p}^{-1}\left(t\right)a_{p}\right\Vert _{\infty}\leq1\right\} $
and the controls for the evader are constrained to the set $\mathcal{A}_{e}=\left\{ a_{e}:\left\Vert Q_{e}^{-1}a_{e}\right\Vert _{\infty}\leq1\right\} $.
The infinity norm with diagonal matrix $Q$, scales the control limit
independently in orthogonal directions. $Q_{p}$ is a function of
time since some systems have control bounds that vary with time, and
is needed to model aerodynamic control surfaces on decelerating vehicles.
Both controls are considered symmetric (centered at zero) for this
paper and all simulations.

We represent the capture set, $\Omega$, as an ellipsoid
\begin{equation}
\Omega=\left\{ x:\left\langle x,W^{-1}x\right\rangle \leq1\right\} .\label{eq:Capture set}
\end{equation}
where $W$ is the ellipsoid shape matrix. The elements of $W$ are
selected such that the pursuing vehicle must be within a distance
$r$
\[
\left\Vert \left[\begin{array}{c}
\delta x\\
\delta y
\end{array}\right]\right\Vert \leq r,
\]
and the velocity at intercept is within some large bound $V_{\text{max}}$
(we don't care what the velocity was at capture, just as long as capture
has occurred). This gives
\[
W=\left[\begin{array}{cccc}
r^{2} &  & \cdots & 0\\
 & r^{2} &  & \vdots\\
\vdots &  & V_{\text{max}}^{2}\\
0 & \cdots &  & V_{\text{max}}^{2}
\end{array}\right].
\]

\subsection{Multi-Vehicle Model}

For a multi-vehicle problem with $k$ pursuers against a single evader,
the joint state space with state vector $\chi\in\mathbb{R}^{4\times k}$
can be constructed as follows
\begin{align}
\chi=\left[\begin{array}{c}
\dot{x_{1}}\\
\dot{x_{2}}\\
\vdots\\
\dot{x_{k}}
\end{array}\right] & =\left[\begin{array}{cccc}
A &  & \cdots & 0\\
 & A &  & \vdots\\
\vdots &  & \ddots\\
0 & \cdots &  & A
\end{array}\right]\left[\begin{array}{c}
x_{1}\\
x_{2}\\
\vdots\\
x_{k}
\end{array}\right]\nonumber \\
 & +\left[\begin{array}{cccc}
B_{1} &  & \cdots & 0\\
 & B_{2} &  & \vdots\\
\vdots &  & \ddots\\
0 & \cdots &  & B_{k}
\end{array}\right]\left[\begin{array}{c}
a_{p1}\\
a_{p2}\\
\vdots\\
a_{pk}
\end{array}\right]\label{eq:Joint state space system}\\
 & +\left[\begin{array}{c}
D\\
D\\
\vdots\\
D
\end{array}\right]a_{e}\nonumber \\
\implies\dot{\chi} & =\hat{A}\chi+\hat{B}a_{p}+\hat{D}a_{e}.
\end{align}
Collaborative guidance is induced by noticing that capture can happen
by \emph{any} single vehicle of the $k$ vehicles in the system. The
capture set for the $i$-th vehicle in the joint system $\left(\ref{eq:Joint state space system}\right)$
is denoted as
\[
\Omega_{i}=\left\{ \chi:\left\langle \chi,W_{i}^{-1}\chi\right\rangle \leq1\right\} ,
\]
with the shape matrix defined as the block diagonal matrix with $W$
on the $i$-th block of the matrix, and the $4\times4$ matrix $\Sigma=V_{\text{max}}^{2}I$
occupying all other blocks. This implies that the capture set for
the joint system is
\begin{equation}
\Omega=\cup_{i}\Omega_{i}.\label{eq:Capture set as union}
\end{equation}

\section{Hamilton\textendash Jacobi Equations with Bounded Control\label{sec:Hamiltonian-Jacobi-Equations}}

\subsection{Viscosity Solutions with the Hopf Formula}

To compute optimal guidance, we use the generalized Hopf formula \cite{darbon2016algorithms,hopf1965generalized,lions1986hopf}.
Consider system dynamics represented as
\begin{equation}
\dot{x}\left(t\right)=f\left(u\left(t\right)\right)\label{eq: Basic system}
\end{equation}
where $x\left(t\right)\in\mathbb{R}^{n}$ is the system state and
$u\left(t\right)\in\mathcal{C}\subset\mathbb{R}^{m}$ is the control
input, constrained to lie in the convex admissible control set $\mathcal{C}$.
We consider a cost functional for a given initial time $t$, and terminal
time $T$ 
\begin{equation}
K\left(x,t,u\right)=\int_{t}^{T}L\left(u\left(s\right)\right)ds+J\left(x\left(T\right)\right),\label{eq: Cost Function}
\end{equation}
where $x\left(T\right)$ is the solution of $\left(\ref{eq: Basic system}\right)$
at terminal time, $T$. We assume that the terminal cost function
$J:\mathbb{R}^{n}\rightarrow\mathbb{R}$ is convex. The function $L:\mathbb{R}^{n}\rightarrow\mathbb{R}\cup\left\{ +\infty\right\} $
is the running cost, and is assumed proper, lower semicontinuous,
convex, and 1-coercive. The value function $v:\mathbb{R}^{n}\times(-\infty,T]\rightarrow\mathbb{R}$
is defined as the minimum cost, $K$, among all admissible controls
for a given state $x$, and time $t\leq T$ with
\begin{equation}
v\left(x,t\right)=\underset{u\in\mathcal{C}}{\text{inf}}\,K\left(x,t,u\right).\label{eq: Value function}
\end{equation}
The value function in $\left(\ref{eq: Value function}\right)$ satisfies
the dynamic programming principle \cite{bryson1975applied,evans10}
and also satisfies the following initial value Hamilton-Jacobi (HJ)
equation by defining the function $\varphi:\mathbb{R}^{n}\times\rightarrow\mathbb{R}$
as $\varphi\left(x,t\right)=v\left(x,T-t\right)$, with $\varphi$
being the viscosity solution of
\begin{equation}
\begin{cases}
\frac{\partial\varphi}{\partial t}\left(x,t\right)+H\left(t,\nabla_{x}\varphi\left(x,t\right)\right)=0 & \text{in}\,\mathbb{R}^{n}\times\left(0,+\infty\right),\\
\varphi\left(x,0\right)=J\left(x\right) & \forall x\in\mathbb{R}^{n},
\end{cases}\label{eq:Initial value HJ PDE}
\end{equation}
where the Hamiltonian $H:\mathbb{R}^{n}\rightarrow\mathbb{R}\cup\left\{ +\infty\right\} $
is defined by
\begin{equation}
H\left(p\right)=\underset{c\in\mathbb{R}^{m}}{\text{sup}}\left\{ \left\langle -f\left(c\right),p\right\rangle -L\left(c\right)\right\} .\label{eq: Basic Hamiltonian definition}
\end{equation}
To apply the constraint that the control must bounded, we introduce
the following running cost $L=\mathcal{I}_{\mathcal{C}}$, where
\[
\mathcal{I}_{\mathcal{C}}=\begin{cases}
0 & \text{if}\,c\in\mathcal{C}\\
+\infty & \text{otherwise,}
\end{cases}
\]
is the indicator function for the set $\mathcal{C}$. This induces
a time-optimal control formulation and reduces the Hamiltonian to
\[
H\left(p\right)=\underset{c\in\mathcal{C}}{\text{max}}\left\langle -f\left(c\right),p\right\rangle .
\]
Solving the HJ equation $\left(\ref{eq:Initial value HJ PDE}\right)$
describes how the value function evolves with time at any point in
the state space and from this, optimal control policies can be found.

It was shown in \cite{darbon2016algorithms} that an exact, point-wise
viscosity solution to $\left(\ref{eq:Initial value HJ PDE}\right)$
can be found using the Hopf formula \cite{hopf1965generalized}. The
value function can be found with the Hopf formula
\begin{equation}
\varphi\left(x,t\right)=-\underset{p\in\mathbb{R}^{n}}{\text{min}}\left\{ J^{\star}\left(p\right)+tH\left(p\right)-\left\langle x,p\right\rangle \right\} ,\label{eq: Basic Hopf formula}
\end{equation}
where the Fenchel-Legendre transform $g^{\star}:\mathbb{R}^{n}\rightarrow\mathbb{R}\cup\left\{ +\infty\right\} $
of a convex, proper, lower semicontinuous function $g:\mathbb{R}^{n}\rightarrow\mathbb{R}\cup\left\{ +\infty\right\} $
is defined by \cite{ekeland1999convex} 
\begin{equation}
g^{\star}\left(p\right)=\underset{x\in\mathbb{R}^{n}}{\text{sup}}\left\{ \left\langle p,x\right\rangle -g\left(x\right)\right\} .\label{eq: Fenchel transform}
\end{equation}
Following the basic definition of the Fenchel-Legendre transform,
$\left(\ref{eq: Basic Hopf formula}\right)$ can be written \cite{lions1986hopf}
as
\[
\varphi\left(x,t\right)=\left(J^{\star}+tH\right)^{\star}\left(x\right).
\]
This shows that value function is itself a Fenchel-Legendre transform.
It follows from a well known property of the Fenchel-Legendre transform
\cite{darbon2015convex} that the unique minimizer of $\left(\ref{eq: Basic Hopf formula}\right)$
is the gradient of the value function
\[
\nabla_{x}\varphi\left(x,t\right)=\text{arg}\,\underset{p\in\mathbb{R}^{n}}{\text{min}}\left\{ J^{\star}\left(p\right)+tH\left(p\right)-\left\langle x,p\right\rangle \right\} ,
\]
provided the gradient exists. So by solving for the value function
using $\left(\ref{eq: Basic Hopf formula}\right)$, we automatically
solve for the gradient.

\subsection{General Linear Models}

Now consider the following linear state space model
\begin{equation}
\dot{x}\left(t\right)=Ax\left(t\right)+B\left(t\right)u\left(t\right),\label{eq:general linear system}
\end{equation}
with $A\in\mathbb{R}^{n\times n}$, $B\in\mathbb{R}^{n\times m}$,
state vector $x\in\mathbb{R}^{n}$, and control input $u\in\mathbb{\mathcal{C}\subset R}^{m}$.
We can make a change of variables
\begin{equation}
z\left(t\right)=e^{-tA}x\left(t\right),\label{eq:change of varibles}
\end{equation}
which results in the following system
\begin{equation}
\dot{z}\left(t\right)=e^{-tA}B\left(t\right)u\left(t\right),\label{eq:z transformed system}
\end{equation}
with terminal cost function now defined in $z$ with $\varphi\left(z,0\right)=J_{z}\left(z,0\right)=J_{x}\left(e^{TA}z\right)$,
which depends on terminal time, $T$. Notice that the system is of
the form presented in $\left(\ref{eq: Basic system}\right)$, with
the exception that the system is now time-varying. It was shown in
\cite[Section 5.3.2, p. 215]{kurzhanski2014dynamics} that the Hopf
formula in $\left(\ref{eq: Basic Hopf formula}\right)$ can be generalized
for a time-varying Hamiltonian to find the value function of the system
in $\left(\ref{eq:z transformed system}\right)$ with
\begin{equation}
\varphi\left(z,t\right)=-\underset{p\in\mathbb{R}^{n}}{\text{min}}\left\{ J_{z}^{\star}\left(p,t\right)+\int_{0}^{t}H\left(p,s\right)ds-\left\langle z,p\right\rangle \right\} ,\label{eq:generalized hopf formula}
\end{equation}
with the time-varying Hamiltonian defined as
\[
H\left(p,t\right)=\underset{c\in\mathcal{C}}{\text{max}}\left\langle e^{-\left(T-t\right)A}B\left(T-t\right)c,p\right\rangle .
\]
The change of variable to $\left(T-t\right)$ is required for time
since the problem was converted to an initial value formulation from
a terminal value formulation in $\left(\ref{eq:Initial value HJ PDE}\right)$. 

\subsection{Linear Differential Games}

Now consider the system
\begin{equation}
\dot{x}\left(t\right)=Ax\left(t\right)+B\left(t\right)u(t)+D\left(t\right)w\left(t\right),\label{eq:linear differential game}
\end{equation}
with $D\left(t\right)\in\mathbb{R}^{n\times\ell}$, which is equal
to $\left(\ref{eq:general linear system}\right)$ with an extra term,
$D\left(t\right)w\left(t\right)$, added. We assume that the additional
control input $w\left(t\right)$ is adversarial and bounded by $w\left(t\right)\in\mathcal{D}\subset\mathbb{R}^{\ell}$.
The cost functional becomes
\begin{equation}
G\left(x,t,u,w\right)=\int_{t}^{T}L\left(u\left(t\right),w\left(t\right)\right)dt+J\left(x\left(T\right)\right),\label{eq:cost functional diff games}
\end{equation}
where $x\left(T\right)$ is the solution of $\left(\ref{eq:linear differential game}\right)$
at terminal time, $T$. We assume that the goal of the adversarial
control input $w\left(t\right)$ is to \emph{increase} the cost functional
$\left(\ref{eq:cost functional diff games}\right)$, in direct contradiction
with the input $u\left(t\right)$, which we are designing in an attempt
to minimize the cost. This system forms a differential game \cite{isaacs1999differential},
and has a corresponding lower value function
\[
V\left(x,t\right)=\underset{u\in\mathcal{C}}{\text{inf}}\,\underset{w\in\mathcal{D}}{\text{sup}}\,G\left(x,t,u,w\right),
\]
and upper value function
\[
U\left(x,t\right)=\underset{w\in\mathcal{D}}{\text{sup}}\,\underset{u\in\mathcal{C}}{\text{inf}}\,G\left(x,t,u,w\right).
\]
As derived in \cite{evans1983differential}, the upper and lower value
functions are viscosity solutions of possibly non convex HJ equation.
We can define the following upper and lower Hamiltonians as
\begin{align*}
H^{+}\left(p,t\right) & =\underset{c\in\mathbb{R}^{m}}{\text{sup}}\,\underset{d\in\mathbb{R}^{\ell}}{\text{inf}}\left\{ \left\langle -f\left(t,c,d\right),p\right\rangle -L\left(c,d\right)\right\} ,\\
H^{-}\left(p,t\right) & =\underset{d\in\mathbb{R}^{\ell}}{\text{inf}}\,\underset{c\in\mathbb{R}^{m}}{\text{sup}}\left\{ \left\langle -f\left(t,c,d\right),p\right\rangle -L\left(c,d\right)\right\} .
\end{align*}
The running cost becomes
\[
L\left(u,w\right)=\mathcal{I}_{\mathcal{C}}\left(u\right)-\mathcal{I}_{\mathcal{D}}\left(w\right),
\]
where $\mathcal{I}_{\mathcal{D}}$ is the indicator function of the
convex set $\mathcal{D}$. If the Hamiltonians $H^{+}$ and $H^{-}$
coincide, then from \cite{evans1983differential}
\[
H^{+}\left(p,t\right)=H^{-}\left(p,t\right)=H^{\pm}\left(p,t\right)\implies U\left(x,t\right)=V\left(x,t\right).
\]
We can apply the same change of variables from $\left(\ref{eq:change of varibles}\right)$
to get
\begin{equation}
\dot{z}\left(t\right)=e^{-tA}B\left(t\right)u\left(t\right)+e^{-tA}D\left(t\right)w\left(t\right),\label{eq:transformed game system}
\end{equation}
and then we can find a candidate solution of the value function $\varphi\left(z,t\right)=U\left(z,t\right)=V\left(z,t\right)$
with the generalized Hopf formula
\[
\varphi\left(z,t\right)=-\underset{p\in\mathbb{R}^{n}}{\text{min}}\left\{ J_{z}^{\star}\left(p,t\right)+\int_{0}^{t}H^{\pm}\left(p,s\right)ds-\left\langle z,p\right\rangle \right\} ,
\]
with the time-varying, non convex Hamiltonian given by
\begin{align}
H^{\pm}\left(p,t\right) & =\underset{c\in\mathcal{C}}{\text{max}}\left\langle e^{-\left(T-t\right)A}B\left(T-t\right)c,p\right\rangle \nonumber \\
 & -\underset{d\in\mathcal{D}}{\text{max}}\left\langle e^{-\left(T-t\right)A}D\left(T-t\right)d,p\right\rangle .\label{eq:Hopf for game}
\end{align}
In general, if $H^{+}\left(p,t\right)\neq H^{-}\left(p,t\right)$,
then the Hopf formula in $\left(\ref{eq:Hopf for game}\right)$ does
not hold. 

\section{Time-Optimal Control with the Hopf Formula\label{sec:Collaborative-Guidance}}

Following the methods presented above in $\left(\ref{eq:transformed game system}\right)$,
we have the transformed system $\left(\ref{eq:Joint state space system}\right)$
as
\[
\dot{z}\left(t\right)=e^{-t\hat{A}}\hat{B}a_{p}\left(t\right)+e^{-t\hat{A}}\hat{D}a_{e}\left(t\right),
\]
and the Hamiltonian is the dual norm of the control set
\begin{align}
H\left(p,t\right) & =\left\Vert Q_{p}\left(T-t\right)\hat{B}^{\dagger}e^{-\left(T-t\right)\hat{A}^{\dagger}}p\right\Vert _{1}\label{eq:game hamiltonian}\\
 & -\left\Vert Q_{e}\hat{D}^{\dagger}e^{-\left(T-t\right)\hat{A}^{\dagger}}p\right\Vert _{1},\nonumber 
\end{align}
where we denote by $\left\Vert \left(\cdot\right)\right\Vert _{1}$
the 1-norm. We choose a convex terminal cost function $J\left(z,0\right)$
such that
\begin{equation}
\begin{cases}
J\left(z,0\right)<0 & \text{for any}\,z\in\text{int}\,\Omega,\\
J\left(z,0\right)>0 & \text{for any}\,z\in\left(\mathbb{R}^{n}\setminus\Omega\right),\\
J\left(z,0\right)=0 & \text{for any}\,z\in\left(\Omega\setminus\text{int}\,\Omega\right),
\end{cases}\label{eq:Terminal Set Def}
\end{equation}
 where $\text{int}\,\Omega$ denotes the interior of $\Omega$. The
intuition behind defining the terminal cost function this way is simple.
If the value function $\varphi\left(z_{0},T\right)<0$ for some $z_{0}$
and $T$, then there exists a control $u\left(t\right)$ that drives
the state from the initial condition at $z_{0}$, to the final state,
$z\left(T\right)$ inside the set $\Omega$. The smallest value of
time $T$, such that $\varphi\left(z_{0},T\right)=0$ is the minimum
time to reach the set $\Omega$, starting at state $z_{0}$. The control
associated with the minimum time to reach is the time-optimal control.
The ellipsoid terminal set defined in $\left(\ref{eq:Capture set}\right)$
results in a quadratic terminal cost function
\[
J_{x}\left(x\right)=\left\langle x,W^{-1}x\right\rangle -1,
\]
After variable substitution the cost function becomes
\[
J_{z}\left(z\right)=\left\langle z,V\left(0\right)z\right\rangle -1,
\]
with $V\left(t\right)=e^{\left(T-t\right)\hat{A}^{\dagger}}W^{-1}e^{\left(T-t\right)\hat{A}}$.
Following the property that the Fenchel-Legendre transform of a norm
function is the dual norm \cite{boyd2004convex}, we have
\[
J_{z}^{\star}\left(p,t\right)=1+\frac{1}{4}\left\langle p,V\left(0\right)^{-1}p\right\rangle .
\]

The generalized Hopf formula requires the integration of the Hamiltonian
which is approximated by Riemann sum quadrature \cite{anton1999calculus}
with step size $h$
\begin{align*}
\int_{0}^{t}H\left(p,t\right)ds & \approx h\sum_{s_{k}\in\mathcal{S}}H\left(p,s_{k}\right),
\end{align*}
where $\mathcal{S}$ denotes the set of discrete time samples. Rectangular
quadrature with fixed step size $h$ was used to pre-compute the time
samples $s_{k}$ from time $0$ to $T$, which requires only a simple
sum at run time to evaluate the integral. We can approximate the matrix
exponential terms efficiently at fixed time intervals, with bounded
error, using \cite{al2011computing}.

To solve the Hopf formula in $\left(\ref{eq:generalized hopf formula}\right)$,
we are performing an unconstrained minimization problem where the
objective function is non-smooth. Non-smooth unconstrained minimization
problems can be solved in a variety of ways. However, because we can
explicitly derive the gradient and Hessian, this directs the use of
a relaxed Newton's method \cite{chapra1998numerical}. We chose for
the initial guess of Newton's method $p_{0}=\frac{V\left(0\right)z}{2}$,
the minimum of the Hopf objective without the Hamiltonian integral.
The initial step size is 1 (full Newton), and is halved whenever the
function value increases during an iteration (without updating the
search direction). The minimization is terminated when the norm of
the change in iterations is small. Most importantly for efficient
implementation, the gradient and Hessian (ignoring discontinuities),
denoted as $\nabla_{p}$ and $\mathcal{H}_{p}$, respectively, for
the minimization can be found directly. The gradient is
\begin{align*}
\nabla_{p}\varphi\left(z,t\right) & =\frac{V\left(0\right)^{-1}p}{2}-z\\
 & +h\sum_{s_{k}\in\mathcal{S}}\Big(Q_{p}\left(s_{k}\right)E_{p}\left(s_{k}\right)\text{sgn}\left(E_{p}\left(s_{k}\right)^{\dagger}p\right)\\
 & -Q_{e}\left(s_{k}\right)E_{e}\left(s_{k}\right)\text{sgn}\left(E_{e}\left(s_{k}\right)^{\dagger}p\right)\Big),
\end{align*}
with $E_{p}\left(t\right)=e^{-\left(T-t\right)\hat{A}}\hat{B}$, and
$E_{e}\left(t\right)=e^{-\left(T-t\right)\hat{A}}\hat{D}$. Additionally
the Hessian is
\[
\mathcal{H}_{p}\left(\varphi\left(z,t\right)\right)=\frac{V\left(0\right)^{-1}}{2}.
\]

To find the optimal control to the desired convex terminal set $\Omega$,
we proceed by solving for the $T^{*}$, the minimum time to reach
the boundary of the set $\Omega$. This is solved numerically with
\[
T^{*}=\text{arg}\,\underset{t<T}{\text{min}}\,\varphi\left(z_{0},t\right).
\]
If the minimum time to reach $T^{*}$ is greater than total available
time $T$, then the set $\Omega$ is not reachable in time $T$ .
The optimal control can then be found from the following relation

\[
\nabla_{p}H\left(\nabla_{z}\varphi\left(z_{0},T^{*}\right),T^{*}\right)=e^{-t\hat{A}}\hat{B}\left(t\right)a_{p}^{*}+e^{-t\hat{A}}\hat{D}\left(t\right)a_{e}^{*}.
\]

To induce collaborative guidance we proceed to solve for the joint
terminal set in $\left(\ref{eq:Capture set as union}\right)$. Let
$J_{i}$ represent terminal cost function of vehicle $i$ with shape
matrix $W_{i}$, then the terminal cost function of the collaborative
system is 
\begin{equation}
J\left(z,t\right)=\underset{i=1,\ldots,k}{\text{min}}J_{i}\left(z,t\right).\label{eq:terminal as union}
\end{equation}
It was shown in \cite{darbon2016algorithms} that max/min-plus algebra
\cite{akian2006max,fleming1997deterministic,mceneaney2006max} can
be used to generalize the Hopf formula to solve for non-convex initial
data that can be formed as the union of convex sets, such as the terminal
cost considered in $\left(\ref{eq:Capture set as union}\right)$.
This is true provided that the Hamiltonian is convex. In general,
the Hamiltonian of the differential game given in $\left(\ref{eq:game hamiltonian}\right)$
is non-convex. But consider the case where $Q_{e}\leq Q_{p}$ and
the system is constrained to the form in $\left(\ref{eq:linearized system}\right)$
and $\left(\ref{eq:Joint state space system}\right)$, then $\left(\ref{eq:game hamiltonian}\right)$
is convex and max/min-plus algebra holds. To find the value function
with the terminal set given by $\left(\ref{eq:terminal as union}\right)$,
we solve the $k$ initial value problems of the form
\begin{equation}
\begin{cases}
\frac{\partial\phi_{i}}{\partial t}\left(z,t\right)+H\left(t,\nabla_{z}\phi_{i}\left(z,t\right)\right)=0 & \text{in}\,\mathbb{R}^{n}\times\left(0,+\infty\right),\\
\phi_{i}\left(z,0\right)=J_{i}\left(z\right) & \forall z\in\mathbb{R}^{n},
\end{cases}\label{eq:max/min-plus algebra}
\end{equation}
and take the pointwise minimum over the $k$ solutions $\phi_{i}\left(z,t\right)$,
each of which has convex initial data, with
\[
\varphi\left(z,t\right)=\underset{i=1,\ldots,k}{\text{min}}\phi_{i}\left(z,t\right).
\]
Each $\phi_{i}\left(z,t\right)$ in $\left(\ref{eq:max/min-plus algebra}\right)$
are independent of each other, and can be computed in parallel. In
the case where the $\left(\ref{eq:game hamiltonian}\right)$ is non-convex,
then the pointwise minimum is only an upper bound of the true value
function; see \cite{m2016idempotent} for more details. 

\section{Results\label{sec:Results}}

\begin{figure}
\begin{centering}
\includegraphics[width=8cm]{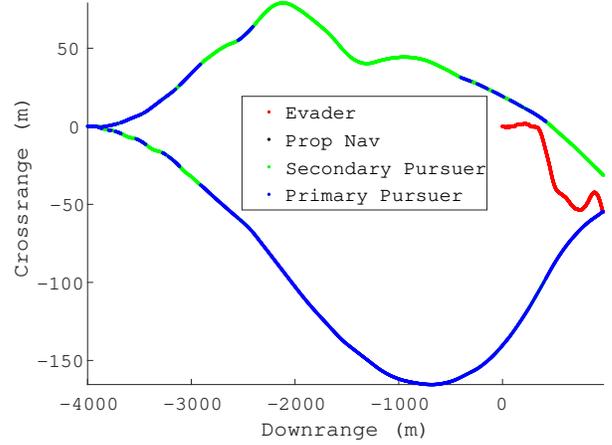}
\par\end{centering}
\caption{The trajectory of Example 1, a tail chase scenario. The red is the
trajectory of the evader. The pursuing vehicles are shown in green
and blue. Blue indicates at that time, it was the pointwise minimum
of the $k$ initial value problems in $\left(\ref{eq:max/min-plus algebra}\right)$,
while green was not. \label{fig:example 1}}
\end{figure}

The above control solution has been integrated into a closed loop
2-on-1 pursuit-evasion 3 degree of freedom (3DOF) simulation using
MATLAB R2016a and Simulink at 120Hz with Euler integration. This included
using a third order autopilot for each pursuer, and using the gradient
of the value function to find optimal evader control. Preliminary
results solved for the optimal control on average $40-83\,ms$ on
a 3 GHz Intel Core i7 950.

As a post-process, the evader\textquoteright s inertial state is found
by solving the modified Dubin's car initial-value problem $\left(\ref{eq:Original non-linear system}\right)$
relative to a fixed origin with zero initial conditions and known
inputs. Adding the evader\textquoteright s inertial state to the vehicle\textquoteright s
relative state and correcting for the induced rotational motion provides
the vehicle\textquoteright s inertial state.

The first example uses a simple geometry in the tail-chase scenario
and the engagement trajectory is shown in Figure \ref{fig:example 1}.
The capture radius is $r=3\,m$, evader control is limited to $\left\Vert a_{e}\right\Vert \leq10\,m/s^{2}$,
and both pursuers have control bounds that decrease in time with
\[
\left\Vert a_{p}\right\Vert \leq\frac{\left(t-40\right)^{2}}{40}\,m/s^{2},
\]
when $0\leq t\leq40$, and $0$ otherwise. The evader is assumed to
travel at speed $V_{e}=50\,m/s$ and the pursuers at speed $255.225\,m/s$
($0.75$ Mach). Both pursuing vehicles, initially launched at $4000\,m$
from the evader, are simultaneously traveling directly at the evader.
Notice that both pursuers separate as to surround and contain the
evader. The miss distance was $0.879\,m<r=3\,m$ and time to intercept
was $19.533$ seconds. In this example, $Q_{e}\leq Q_{p}$ and the
Hamiltonian remained convex for the duration of the simulated engagement. 

\begin{figure}
\begin{centering}
\includegraphics[width=8cm]{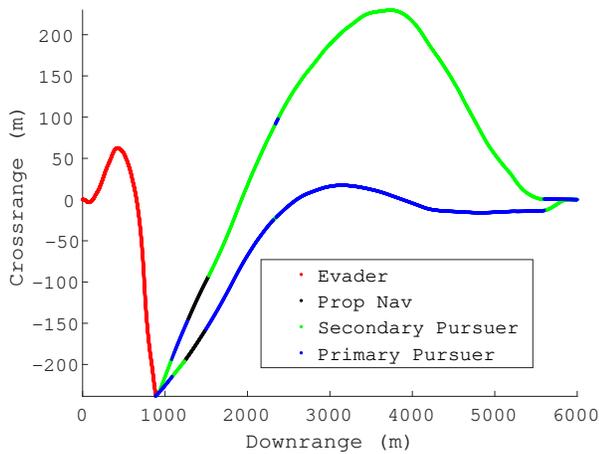}
\par\end{centering}
\caption{The trajectory of Example 2, a head-on scenario. The red is the trajectory
of the evader. The pursuing vehicles are shown in green and blue.
Blue indicates at that time, it was the pointwise minimum of the $k$
initial value problems in $\left(\ref{eq:max/min-plus algebra}\right)$,
while green was not. Black is when $\Omega$ was considered not reachable
due to high linearization error and proportional navigation was used.\label{fig:Example 2}}
\end{figure}

The second example utilized a similar engagement, but with head-on
aspect configuration. The parameters are the same as example 1, but
a $6000\,m$ initial separation. In this case, the initial conditions
are such that during simulation, the linearization error in $\left(\ref{eq:linearized system}\right)$
is large. When this occurs, the solution of the zero level set time
maybe higher than available flight time $T$. This indicates the set
$\Omega$ is not reachable (due to the linearization error) and in
our simulations reverts to proportional navigation (PN) until the
set $\Omega$ is considered reachable. This can easily be countered
by increasing the control bound of the evader to account for linearization
error. Additionally, the convexity assumption of $Q_{e}\leq Q_{p}$
is violated in this example, but only for the last $0.158$ seconds,
or about $0.78\%$ of the engagement. With both vehicles launched
simultaneously, intercept still occurred, with a miss distance of
$2.34\,m$. Time to intercept was $20.158$ seconds and the flyout
paths are given in Figure \ref{fig:Example 2}. 

\section{Conclusions and Future Work}

The generalized Hopf formula provide new capabilities for solving
high-dimensional optimal control and differential games, such as the
pursuit-evasion guidance presented here. Additionally, the above work
can be used for evasion strategies that could be of interest for collision
avoidance problems. Future work will focus on extending the generalized
Hopf formula for certain classes of non-linear systems, such as feedback
linearizable systems \cite{slotine1991applied}, and apply splitting
algorithms \cite{darbon2016algorithms,goldstein2009split,chambolle2011first}
for efficient optimization when the gradient and Hessian is not explicitly
known.

\section*{Acknowledgments}

The authors would like to thank the anonymous reviewers. Their comments
and suggestions greatly improved the accuracy and clarity of this
paper. 

\bibliographystyle{plain}
\bibliography{CDC2017}

\begin{thebibliography}{10}

\bibitem{akian2006max}
M.~Akian, R.~Bapat, and S.~Gaubert.
\newblock Max-plus algebra.
\newblock {\em Handbook of Linear Algebra (Discrete Mathematics and its
  Applications)}, 39:10--14, 2006.

\bibitem{al2011computing}
A.~H. Al-Mohy and N.~J. Higham.
\newblock Computing the action of the matrix exponential, with an application
  to exponential integrators.
\newblock {\em SIAM Journal on Scientific Computing}, 33(2):488--511, 2011.

\bibitem{anton1999calculus}
H.~Anton, S.~Davis, and I.~Bivens.
\newblock {\em Calculus: A New Horizon}.
\newblock Wiley New York, 1999.

\bibitem{bellman1957dynamic}
R.~E. Bellman.
\newblock {\em Dynamic Programming}, volume~1.
\newblock Princeton University Press, 1957.

\bibitem{bellman2015adaptive}
R.~E. Bellman.
\newblock {\em Adaptive Control Processes: A Guided Tour}.
\newblock Princeton University Press, 2015.

\bibitem{boyd2004convex}
S.~Boyd and L.~Vandenberghe.
\newblock {\em Convex Optimization}.
\newblock Cambridge University Press, 2004.

\bibitem{bryson1975applied}
A.~R. Bryson and Y.-C. Ho.
\newblock {\em Applied Optimal Control: Optimization, Estimation and Control}.
\newblock CRC Press, 1975.

\bibitem{chambolle2011first}
A.~Chambolle and T.~Pock.
\newblock A first-order primal-dual algorithm for convex problems with
  applications to imaging.
\newblock {\em Journal of Mathematical Imaging and Vision}, 40(1):120--145,
  2011.

\bibitem{chapra1998numerical}
S.~C. Chapra and R.~P. Canale.
\newblock {\em Numerical Methods for Engineers}, volume~2.
\newblock McGraw-Hill New York, 1998.

\bibitem{darbon2015convex}
J.~Darbon.
\newblock On convex finite-dimensional variational methods in imaging sciences
  and {Hamilton}-{Jacobi} equations.
\newblock {\em SIAM Journal on Imaging Sciences}, 8(4):2268--2293, 2015.

\bibitem{darbon2016algorithms}
J.~Darbon and S.~Osher.
\newblock Algorithms for overcoming the curse of dimensionality for certain
  {Hamilton}-{Jacobi} equations arising in control theory and elsewhere.
\newblock {\em Research in the Mathematical Sciences}, 3(1):19, 2016.

\bibitem{dubins1957curves}
L.~E. Dubins.
\newblock On curves of minimal length with a constraint on average curvature,
  and with prescribed initial and terminal positions and tangents.
\newblock {\em American Journal of Mathematics}, 79(3):497--516, 1957.

\bibitem{ekeland1999convex}
I.~Ekeland and R.~Temam.
\newblock {\em Convex Analysis and Variational Problems}.
\newblock SIAM, 1999.

\bibitem{evans1983differential}
L.~C. Evans and P.~E. Souganidis.
\newblock Differential games and representation formulas for solutions of
  {Hamilton}-{Jacobi}-{Isaacs} equations.
\newblock Technical report, DTIC Document, 1983.

\bibitem{evans10}
Lawrence~C. Evans.
\newblock {\em Partial differential equations}.
\newblock American Mathematical Society, Providence, R.I., 2010.

\bibitem{fleming1997deterministic}
W.~H. Fleming.
\newblock Deterministic nonlinear filtering.
\newblock {\em Annali della Scuola Normale Superiore di Pisa-Classe di
  Scienze}, 25(3-4):435--454, 1997.

\bibitem{goldstein2009split}
T.~Goldstein and S.~Osher.
\newblock The split {Bregman} method for l1-regularized problems.
\newblock {\em SIAM Journal on Imaging Sciences}, 2(2):323--343, 2009.

\bibitem{hopf1965generalized}
E.~Hopf.
\newblock Generalized solutions of non-linear equations of first order.
\newblock {\em Journal of Mathematics and Mechanics}, 14:951--973, 1965.

\bibitem{huang2011differential}
H.~Huang, J.~Ding, W.~Zhang, and C.~J. Tomlin.
\newblock A differential game approach to planning in adversarial scenarios: A
  case study on capture-the-flag.
\newblock In {\em 2011 IEEE International Conference on Robotics and Automation
  (ICRA)}, pages 1451--1456. IEEE, 2011.

\bibitem{isaacs1999differential}
R.~Isaacs.
\newblock {\em Differential Games: A Mathematical Theory with Applications to
  Warfare and Pursuit, Control and Optimization}.
\newblock Courier Corporation, 1999.

\bibitem{kurzhanski2014dynamics}
A.~B. Kurzhanski and P.~Varaiya.
\newblock {\em Dynamics and Control of Trajectory Tubes: Theory and
  Computation}, volume~85.
\newblock Springer, 2014.

\bibitem{lions1986hopf}
P.~L. Lions and J.-C. Rochet.
\newblock Hopf formula and multitime {Hamilton}-{Jacobi} equations.
\newblock {\em Proceedings of the American Mathematical Society}, 96(1):79--84,
  1986.

\bibitem{mceneaney2006max}
W.~M. McEneaney.
\newblock {\em Max-Plus Methods for Nonlinear Control and Estimation}.
\newblock Springer Science \& Business Media, 2006.

\bibitem{m2016idempotent}
W.~M. McEneaney and A.~Pandey.
\newblock An idempotent algorithm for a class of network-disruption games.
\newblock {\em Kybernetika}, 52(5):666--695, 2016.

\bibitem{mitchell2004toolbox}
I.~Mitchell.
\newblock A toolbox of level set methods.
\newblock {\em Dept. Comput. Sci., Univ. British Columbia, Vancouver, BC,
  Canada, http://www. cs. ubc. ca/\~{} mitchell/ToolboxLS/toolboxLS. pdf, Tech.
  Rep. TR-2004-09}, 2004.

\bibitem{mitchell2008flexible}
I.~Mitchell.
\newblock The flexible, extensible and efficient toolbox of level set methods.
\newblock {\em Journal of Scientific Computing}, 35(2):300--329, 2008.

\bibitem{mitchell2005time}
I.~Mitchell, A.~M. Bayen, and C.~J. Tomlin.
\newblock A time-dependent {Hamilton}-{Jacobi} formulation of reachable sets
  for continuous dynamic games.
\newblock {\em IEEE Transactions on Automatic Control}, 50(7):947--957, 2005.

\bibitem{osher2006level}
S.~Osher and R.~Fedkiw.
\newblock {\em Level Set Methods and Dynamic Implicit Surfaces}, volume 153.
\newblock Springer Science \& Business Media, 2006.

\bibitem{palumbo2010basic}
N.~F. Palumbo, R.~A. Blauwkamp, and J.~M Lloyd.
\newblock Basic principles of homing guidance.
\newblock {\em Johns Hopkins APL Technical Digest}, 29(1):25--41, 2010.

\bibitem{palumbo2010modern}
N.~F. Palumbo, R.~A. Blauwkamp, and J.~M Lloyd.
\newblock Modern homing missile guidance theory and techniques.
\newblock {\em Johns Hopkins APL Technical Digest}, 29(1):42--59, 2010.

\bibitem{pan2012pursuit}
S.~Pan, H.~Huang, J.~Ding, W.~Zhang, and C.~J. Tomlin.
\newblock Pursuit, evasion and defense in the plane.
\newblock In {\em American Control Conference (ACC)}, pages 4167--4173. IEEE,
  2012.

\bibitem{slotine1991applied}
J.-J. Slotine and W.~Li.
\newblock {\em Applied Nonlinear Control}, volume 199.
\newblock Prentice-Hall Englewood Cliffs, NJ, 1991.

\bibitem{stipanovic2004decentralized}
D.~M. Stipanovi{\'c}, G.~Inalhan, R.~Teo, and C.~J. Tomlin.
\newblock Decentralized overlapping control of a formation of unmanned aerial
  vehicles.
\newblock {\em Automatica}, 40(8):1285--1296, 2004.

\end{thebibliography}

\end{document}